%% file: NLTG_v5_20181227.tex
\documentclass[12pt]{iopart}

\newcommand\eqref[1]{(\ref{#1})}

\usepackage{iopams}
\usepackage{amssymb}
\usepackage{pifont}
\usepackage{mathrsfs}
\usepackage{amsfonts,amssymb,theorem,graphicx,listings,txfonts,multirow}
\usepackage{algorithm}
\usepackage{algpseudocode}
\usepackage{fancyhdr}
\usepackage{subfigure}
\usepackage{soul,color,xcolor}
\usepackage{epstopdf}
\newtheorem{theorem}{Theorem}
\newtheorem{lemma}{Lemma}

\newtheorem{proposition}[theorem]{Proposition}
\newtheorem{corollary}[theorem]{Corollary}
\usepackage[colorlinks,
            linkcolor=red,
            citecolor=blue,
            urlcolor=red
            ]{hyperref}
\newenvironment{pf}{\vspace*{1pt}\par\noindent{\bf Proof.}}%
                    {\vspace{1em}$\blacksquare$} % square <--> blacksquare

\input Definition_Paper.tex

\begin{document}
%% title
\title[NLTG Priors in Medical Image]{Nonlocal TV-Gaussian (NLTG) prior for Bayesian inverse problems with applications to Limited CT Reconstruction}

\author{Didi Lv$^1$, Qingping Zhou$^1$, Jae Kyu Choi$^{3}$, Jinglai Li$^{4}$, Xiaoqun Zhang$^{1,2,5}$}

\address{
$^1$School of Mathematical Sciences, Shanghai Jiao Tong University, Shanghai 200240, People's Republic of China \\
$^2$Institute of Natural Sciences, Shanghai Jiao Tong University, Shanghai 200240, People's Republic of China \\
$^3$School of Mathematical Sciences, Tongji University, Shanghai 200082, People's Republic of China \\
$^4$Department of Mathematical Sciences, University of Liverpool, Liverpool L69 6ZL, United Kingdom \\
$^5$MOE Key Laboratory of Scientific and Engineering Computing, Shanghai Jiao Tong University, Shanghai 200240, People's Republic of China\\

}

% Email Address
\ead{jinglai.li@liverpool.ac.uk, xqzhang@sjtu.edu.cn}

% abstract
\begin{abstract}
Bayesian inference methods have been widely applied in inverse problems,
{largely due to their ability to characterize the uncertainty associated with the estimation results.}
{In the Bayesian framework} the prior distribution of  the unknown  plays an essential role in the Bayesian inference, {and a good prior distribution
can significantly improve the inference results.}
In this paper, we extend the total~variation-Gaussian (TG) prior  in \cite{Z.Yao2016}, and propose a hybrid prior distribution which combines
the nonlocal total variation regularization and the Gaussian (NLTG) distribution. The advantage of the new prior is two-fold.
The proposed prior models both texture and geometric structures present in images through the NLTV. The Gaussian reference measure also provides a flexibility of incorporating structure information  from a reference image. Some theoretical properties are established for the NLTG prior. The  proposed prior is applied to limited-angle tomography reconstruction problem with difficulties of severe data missing.  We compute both MAP and CM estimates through two efficient methods and the numerical experiments validate the advantages and feasibility of the proposed NLTG prior.
\end{abstract}

\maketitle

\section{Introduction}\label{sec:intro}
{Bayesian inference methods~\cite{gelman2014bayesian,kaipio2006statistical} have become a popular tool to solve inverse problems. Such popularity
is largely due to its ability to quantify the solution uncertainties. A typical Bayesian treatment consists of assigning a prior distribution to the unknown parameters and then update the distribution based on the observed data, yielding the posterior distribution. Recently considerable attentions
have been paid to the studies of infinite dimensional Bayesian inverse problems, where the unknowns are functions of space or time, for example,
images.
In particular, a rigorous function space Bayesian inference framework for inverse problems was
developed in \cite{Bayesianperspective}. It should be clear that, as most practical inverse problems are highly ill-posed,
the performance of the Bayesian inference depends critically on the choice of the prior distribution,
and prior modeling plays an essential role in the Bayesian inference method.
For  infinite dimensional problems, the Gaussian measures are arguably the most popular choice of prior distributions,
as it has many theoretical and computational advantages in the infinite dimensional setting~\cite{Bayesianperspective}.
}

{However, in many practical problems, such as medical image reconstruction,
the functions or images that one wants to recover are often subject to sharp jumps or discontinuities.
 The Gaussian prior distributions are typically not suitable for modeling such functions \cite{yao2016tv}.
To this end several non-Gaussian priors have been proposed to model such images, e.g., \cite{vollmer2015dimension}.
 Since these prior distributions differ significantly from Gaussian, many sampling schemes based on
 the Gaussian prior can not be used directly.
 To address the issue, a hybrid prior was proposed in \cite{Z.Yao2016}.
 The hybrid prior is motivated by the total variation (TV) regularization~\cite{rudin1992nonlinear}  in the deterministic setting;
 however, it has been proven in \cite{lassas2004can} that the TV based prior does not converge to a well-defined infinite-dimensional measure as the discretization dimension increases.
The hybrid prior is a combination of the TV term and the Gaussian distribution: it uses a TV term to capture the sharp jumps in the functions and  a Gaussian distribution as a reference measure to make sure that the resulting prior does
 converge to a well-defined probability measure in the function space in the infinite dimensional limit.  }

Nonlocal methods are  another types of popular regularization methods for imaging inverse problems. They are originally proposed for natural image processing to restore repetitive patterns and textures, for example a heuristic copy-paste technique was firstly proposed for  texture synthesis in \cite{Efros99},  a more systematic nonlocal means filter was proposed in ~\cite{A.Buades2005} and a nonlocal variational framework  was established in \cite{G.Gilboa2008}. The main idea of the nonlocal methods is to utilize the similarities present in an image as a weight for restoring, smoothing or regularization. As an extension of TV, nonlocal TV (NLTV) regularization method is  among those popular variational regularization tools due to its flexibility of recovering both texture and geometry patterns for diverse imaging inverse problems, see \cite{X.Zhang2010,X.Zhang2010a,peyre2011nonlocal,liu2015mo} for the applications. It has been demonstrated that in many practical problems,
the NLTV method has better performance than the standard TV, especially for recovering textures and structures of images.
 More definitions and details are present in Section \ref{sec2.1}.

{Inspired by the success of nonlocal regularization, we propose to improve the hybrid TV-Gaussian (TG) prior proposed in \cite{Z.Yao2016}
by in-cooperating the nonlocal methods. The idea is rather straightforward: we shall replace the TV term in the hybrid prior with a NLTV term,
and theoretically  we are able to prove that the resulting new hybrid prior can also lead to well defined posterior distribution in the infinite dimensional setting.
Moreover, the new hybrid prior has the following advantages:
first the NLTV term can better recover textures and structures of images, especially for highly ill-posed or severely data-missing inverse problems;
secondly, we have extended the function space to a larger function space compared to $H^1$ considered in TG prior, for dealing with larger images class;
finally the Gaussian measure provides some freedom to incorporate other prior information through the covariance matrix, such as structures from a reference image.  To demonstrate the effectiveness of the hybrid NLTV-Gaussian (NLTG) prior,  we apply it to the limited tomography problems, where only limited projection data are available. In particular, we  consider the  two common types of point estimation in the Bayesian framework: maximum a posterior (MAP) and conditional mean (CM) with the NLTG prior. The MAP estimate consists of solving an optimization problem, while the computing of  CM  involves evaluating a high-dimensional integration problem \cite{kaipio2006statistical,lucka2012hierarchical}.} %In this paper, we will present the methods of computing the two types of estimates of NLTG prior for CT reconstruction problem.

{The remainder of this paper is structured as follows: Section 2 describes the proposed  NLTG priors construction on the separable Hilbert space and presents the related theoretical properties. In section 3, we give a simple introduction on the limited tomography and solve the reconstruction problem by applying the proposed NLTG prior. The numerical results via  both MAP  and CM estimates are shown in section 4. Finally, we draw our conclusions in the last section.}

\section{The NLTG priors}

\subsection{The Bayesian framework and the TG prior}

We first give a brief introduction to the basic setup of the Bayesian inference methods for inverse problems.  We consider the forward model of the following form:
\begin{equation}
  y = F(u) + n,
\end{equation}
where $u$ is an unknown function (in this work we shall restrict ourselves to the situation where $u$ is a real-valued function defined in $\R^2$, i.e., an image), $y\in\mathbb{R}^m$ is the measured data and $n$ is a $m$-dimensional zero mean Gaussian noise with covariance matrix $\Sigma_0$.
Our goal here is to estimate the unknown function $u$ from the measured data $y$.

First we assume that the unknown $u$ lives in a Hilbert space of functions, say $\mX$.
We then choose a probabilistic measure defined on $\mX$, denoted by $\mu_{pr}$, to be  the prior measure of $u$.
The posterior measure of $u$, denoted as $\mu^y$, is represented as the Radon-Nikodym (R-N) derivative with respect to $\mu_{pr}$:
\begin{equation}\label{PostMeasure}
\f{\rd\mu^y}{\rd\mu_{pr}}\propto\exp\left(-\Phi(u)\right),
\end{equation}
where $\Phi(u) = \frac12\|F(u)-y\|_{\Sigma_0}^2 =\frac12 \langle\Sigma_0^{-\frac12}(F(u)-y),\Sigma_0^{-\frac12}(F(u)-y)\rangle$
 is  the data fidelity term in deterministic inverse problem.
 In the Bayesian framework, the posterior distribution depends on the information from data and the prior knowledge represented by the prior distribution,
 and so the choice of the prior distribution plays an essential role in the Bayesian method.

As is discussed in Section~\ref{sec:intro}, probably the most popular prior in the infinite dimensional setting is the Gaussian measure.
That is, we choose $\mu_{pr}=\mu_0=N(0,C_0)$  a Gaussian measure defined on $\mX$ with zero mean and covariance operator $C_0$. Note that $V_0$ is symmetric positive operator of trace class \cite{li2015note}.
%The Cameron-Martin space of $\mu_0$ \cite{da2006introduction} is defined as the range of $V_0^{\f{1}{2}}$
%\begin{equation*}
%  \mK = \{u = V_0^{\f{1}{2}}x:x\in\mX\}\subset\mX.
%\end{equation*}
%Then $\mK$ is a Hilbert space equipped with the following inner product
%\begin{equation*}
%  \la\cdot,\ \cdot\ra_\mK=\la V_0^{\f{1}{2}}\cdot,\ V_0^{\f{1}{2}}\cdot\ra_\mX.
%\end{equation*}}
To better model functions with sharp jumps, Ref \cite{Z.Yao2016} proposes the hybrid TG prior in the form of,
\begin{equation}\label{PriorMeasure}
\f{\rd\mu_{pr}}{\rd\mu_0}\propto\exp\left(-R(u)\right),
\end{equation}
where $R(u)$ represents additional prior information (or regularization) on $u$.
In this case, one can writhe the R-N derivative of $\mu^y$ with respect to $\mu_0$:
\begin{equation}\label{RNDerivative}
\f{\rd\mu^y}{\rd\mu_0}\propto\exp\left(-\Phi(u)-R(u)\right),
\end{equation}
which in turn returns to the conventional formulation with Gaussian priors.
Specifically we shall choose the state space $\mX$ to be Sobolev space $H_1$ and $R(u) = \lambda \|u\|_{TV}=\lambda\int\|\nabla u\|_2 \rd\x$ introduced in \cite[Section 2.3]{Z.Yao2016}.

\subsection{Nonlocal Total Variation}\label{NLTV}\label{sec2.1}

Here we provide the formulation of the NLTG prior and we start with a brief introduction to the nonlocal regularization.
For the details, one may consult \cite{G.Gilboa2008} for the variational framework based nonlocal operators, \cite{G.Gilboa2007,Y.Lou2010,X.Zhang2010,X.Zhang2010a} for a short survey on the theory and application of NLTV, and \cite{A.Buades2005,Chung1997,A.Elmoataz2008,S.Kindermann2005,Peyre2008,D.Zhou2005} for more relative surveys.
The nonlocal methods can be described as follows.
Let $\Om\subseteq\R^2$ be a bounded set, and $\x\in\Om$. Given a reference image $f:\Om\to\R$, we define a nonnegative symmetric weight function $\om:\Om\times\Om\to\R$ as follows:
\begin{equation}\label{NonlocalWeight}
\om(\x,\y)=\exp\left\{-\f{\left\la G_a,~\left|f(\x+\cdot)-f(\y+\cdot)\right|^2\right\ra}{h^2}\right\}
\end{equation}
where $G$ is a Gaussian kernel with the standard deviation $a$, $h$ is a filtering parameter {and  $\la\cdot,~\cdot\ra$ is the standard inner product in $L_2(\Omega)$}. Note that in general, $h$ corresponds to the noise level; conventionally we set it to be the standard deviation of the noise \cite{Y.Lou2010}.

Let $u:\Om\to\R$. Using the weight function $\om:\Om\times\Om\to\R$ in (\ref{NonlocalWeight}), we define the nonlocal (NL) gradient $\na_{\om}u:\Om\times\Om\to\R$ as
\begin{equation}\label{NonlocalGradient}
\na_{\om}u(\x,\y)=\left(u(\y)-u(\x)\right)\sqrt{\om(\x,\y)}.
\end{equation}
For a given $p:\Om\times\Om\to\R$, its NL divergence is defined by the standard adjoint relation with the NL gradient operator as follows:
\begin{equation*}
\left\la\na_{\om}u,p\right\ra=-\left\la u,\div_{\om}p\right\ra,
\end{equation*}
which leads to the following explicit formula:
\begin{equation}\label{NonlocalDivergence}
\div_{\om}p(\x)=\int_{\Om}\left(p(\y,\x)-p(\x,\y)\right)\sqrt{\om(\x,\y)}\rd\y.
\end{equation}
Now, we design the following functional based on the nonlocal operators:
\begin{equation}
J_{NLTV}(u)=\int_{\Om}\left\|\na_{\om}u(\x,\cdot)\right\|_{L_2(\Om)}\rd\x=\int_{\Om}\left(\int_{\Om}\left|u(\x)-u(\y)\right|^2\om(\x,\y)\rd\y\right)^{\f{1}{2}}\rd\x.\label{NLTVSeminorm}
\end{equation}
%The discrete version of NLTV can be represented as followed
%\begin{equation}
%  J_{NLTV}(u) = \sum_{\x}\|\sum_{\y}\na_{\om}u(\x,\y)\|_2 = \|\na_{\om}u\|_{1,2}.\label{NLTVSeminormdiscrete}
%\end{equation}

Then we can see that the functional in \eqref{NLTVSeminorm} is analogous to the total variation (TV) seminorm,
and the NLTG prior can be constructed similarly.
To do this, we first need to specify the state space $\mX$, which should be desirably a separable Hilbert space. Recall that $\Om$ is a bounded set in $\R^2$. As convention, we choose $\Om=[0,1]^2$ throughout this paper. (Note, however, that our analysis is valid for any bounded domain with $C^1$ boundary). For a given reference image $f:\Om\to\R$ and a weight function $\om:\Om\times\Om\to\R$ defined as \eqref{NonlocalWeight}, we introduce
\begin{equation}\label{TheSpaceX}
\mX=\left\{u\in L_2(\Om):\na_{\om}u\in L_2(\Om\times\Om)\right\},
\end{equation}
where the NL gradient $\na_{\om}u$ is defined as \eqref{NonlocalGradient}. Obviously, we have $\mX\subseteq L_2(\Om)$. In addition, we present the following lemma which in turn shows that the NLTV functional defined as \eqref{NLTVSeminorm} is well defined on $L_2(\Om)$.
\begin{lemma}\label{Lemma1} Let $f\in L_2(\Om)$ be a given reference image, and let $\om:\Om\times\Om\to\R$ be a weight function defined as \eqref{NonlocalWeight}. For $u\in L_2(\Om)$, we have
\begin{equation}\label{NLGradContinuity}
\left\|\na_{\om}u\right\|_{L_2(\Om\times\Om)}\leq2\left\|u\right\|_{L_2(\Om)}
\end{equation}
and
\begin{equation}\label{NLTVContinuity}
J_{NLTV}(u)\leq2\left\|u\right\|_{L_2(\Om)}
\end{equation}
As a consequence, $\na_{\om}:L_2(\Om)\to L_2(\Om\times\Om)$ is a continuous linear operator.
\end{lemma}

\begin{pf} First we note that, from the definition \eqref{NonlocalWeight} of $\om$, we have
\begin{equation*}
\om(\x,\y)\in(0,1]~~~~~~~\forall\x,\y\in\Om.
\end{equation*}
Hence, using H\"older's inequality (e.g. \cite{Foll1999}), we have
\begin{eqnarray*}
\left\|\na_{\om}u\right\|_{L_2(\Om\times\Om)}^2&=&\int_{\Om}\int_{\Om}\left|\na_{\om}u(\x,\y)\right|^2\rd\y\rd\x\\
&=&\int_{\Om}\int_{\Om}\left|u(\x)-u(\y)\right|^2\om(\x,\y)\rd\y\rd\x\\
&\leq&\int_{\Om}\int_{\Om}\left|u(\x)-u(\y)\right|^2\rd\y\rd\x\\
&\leq&\int_{\Om}\int_{\Om}\left|u(\y)\right|^2\rd\y\rd\x+\int_{\Om}\int_{\Om}\left|u(\x)\right|^2\rd\y\rd\x+2\int_{\Om}\int_{\Om}\left|u(\x)u(\y)\right|\rd\y\rd\x\\
&\leq&2\left\|u\right\|_{L_2(\Om)}^2+2\left\|u\right\|_{L_2(\Om)}^2=4\left\|u\right\|_{L_2(\Om)}^2,
\end{eqnarray*}
where we used the fact that $|\Om|=1$. This concludes \eqref{NLGradContinuity}.

For \eqref{NLTVContinuity}, we first note that for each $\x\in\Om$, the mapping
\begin{equation}
\x\mapsto\left\|\na_{\om}u(\x,\cdot)\right\|_{L_2(\Om)}
\end{equation}
is in $L_2(\Om,\rd\x)$, by \eqref{NLGradContinuity}, where we specified the Lebesgue measure $\rd\x$ for clarity. Then applying the H\"older's inequality again, we have
\begin{equation*}
J_{NLTV}(u)\leq\left(\int_{\Om}\left\|\na_{\om}u(\x,\cdot)\right\|_{L_2(\Om)}^2\rd\x\right)^{\f{1}{2}}=\left\|\na_{\om}u\right\|_{L_2(\Om\times\Om)}\leq2\left\|u\right\|_{L_2(\Om)},
\end{equation*}
and this concludes Lemma \ref{Lemma1}.
\end{pf}

Combining the definition \ref{TheSpaceX} and Lemma \ref{Lemma1}, we can see that $\mX=L_2(\Omega)$. Now by taking
\begin{equation}\label{Regularizer}
R(u)=\lambda J_{NLTV}(u),
\end{equation}
in \eqref{RNDerivative}, we obtain the NLTG prior distribution.

\subsection{Theoretical properties of the NLTG prior}\label{HybridPrior}

In this section, we show that the NLTG prior leads to a well-behaved posterior distribution in $\mX$, where the proofs follow the similar line as \cite{Z.Yao2016}.

{Following \cite{Bayesianperspective}, we assume that the forward operator $F:L_2(\Om)\to\R^m$ satisfies the following conditions:}
\begin{enumerate}
\item[A1.] For every $\ep>0$, there exists $M=M(\ep)>0$ such that
\begin{equation*}
\left\|F(u)\right\|_{\Sigma_0^{-1}}\leq\exp\left(\ep\left\|u\right\|_{L_2(\Om)}^2+M\right)~~~~~\forall u\in L_2(\Om).
\end{equation*}
\item[A2.] For every $\delta>0$, there exists $L=L(\delta)>0$ such that for all $u_1$, $u_2\in L_2(\Om)$ with
\begin{equation*}
\max\left\{\left\|u_1\right\|_{L_2(\Om)},\left\|u_2\right\|_{L_2(\Om)}\right\}\leq\delta,
\end{equation*}
we have
\begin{equation*}
\left\|F(u_1)-F(u_2)\right\|_{\Sigma_0^{-1}}\leq L\left\|u_1-u_2\right\|_{L_2(\Om)}.
\end{equation*}
\end{enumerate}

 We have the following proposition on $\mu_{pr}$ in \eqref{PriorMeasure} .

\begin{proposition}\label{Prop1} Let $R:L_2(\Om)\to\R$ be defined as \eqref{Regularizer}. Then we have the followings:
\begin{enumerate}
\item For all $u\in L_2(\Om)$, we have $R(u)\geq0$.
\item For every $\delta>0$, there exists $M=M(\delta)>0$ such that for all $u\in L_2(\Om)$ with $\left\|u\right\|_{L_2(\Om)}\leq\delta$, we have $R(u)\leq M$.
\item For every $\delta>0$, there exists $L=L(\delta)>0$ such that for all $u_1$, $u_2\in L_2(\Om)$ with
\begin{equation*}
\max\left\{\left\|u_1\right\|_{L_2(\Om)},\left\|u_2\right\|_{L_2(\Om)}\right\}\leq\delta,
\end{equation*}
we have
\begin{equation*}
\left|R(u_1)-R(u_2)\right|\leq L\left\|u_1-u_2\right\|_{L_2(\Om)}.
\end{equation*}
\end{enumerate}
\end{proposition}

\begin{pf} Since (i) is trivial, we focus on (ii) and (iii). For (ii), note that \eqref{NLGradContinuity} of Lemma \ref{Lemma1} gives
\begin{equation*}
J_{NLTV}(u)\leq 2\left\|u\right\|_{L_2(\Om)}.
\end{equation*}
For a given $\delta>0$, we choose $M=2\lambda\delta$. Then whenever $\left\|u\right\|_{L_2(\Om)}\leq\delta$, we have
\begin{equation*}
R(u)=\lambda J_{NLTV}(u)\leq2\lambda\left\|u\right\|_{L_2(\Om)}\leq M,
\end{equation*}
which concludes (ii).
For (iii), first note that for $u_1$, $u_2\in L_2(\Om)$, we have
\begin{eqnarray*}
\left|J_{NLTV}(u_1)-J_{NLTV}(u_2)\right|&=&\left|\int_{\Om}\left\|\na_{\om}u_1(\x,\cdot)\right\|_{L_2(\Om)}\rd\x-\int_{\Om}\left\|\na_{\om}u_2(\x,\cdot)\right\|_{L_2(\Om)}\rd\x\right|\\
&\leq&\int_{\Om}\left\|\na_{\om}\left(u_1-u_2\right)(\x,\cdot)\right\|_{L_2(\Om)}\rd\x\\
&=&J_{NLTV}(u_1-u_2)\leq2\left\|u_1-u_2\right\|_{L_2(\Om)}.
\end{eqnarray*}
where the last inequality follows from \eqref{NLTVContinuity} of Lemma \ref{Lemma1}. Then by letting $L=2\lambda$, we conclude (iii). This completes the proof.
\end{pf}

Theorem~\ref{Th1} states that $\mu^y$ in \eqref{RNDerivative} is a well-defined probability measure on $L_2(\Om)$ and it is Lipschitz continuous in the data $y$. Since the theorem is a direct consequence of the fact that $\Phi+R$ satisfies \cite[Assumptions 2.6.]{Vollmer2015}, we omit the proof.

\begin{theorem}\label{Th1} Assume that $ {F}:L_2(\Om)\to\R^m$ satisfies A1-A2. Let $R:L_2(\Om)\to\R$ be defined as \eqref{Regularizer}. For a given $y\in\R^m$, we define $\mu^y$ as in \eqref{RNDerivative}. Then we have the followings:
\begin{enumerate}
\item $\mu^y$ is a well-defined measure on $L_2(\Om)$.
\item $\mu^y$ is Lipschitz continuous in the data $y$ with respect to the Hellinger distance. More precisely, if $\mu^y$ and $\mu^{y'}$ are two measures corresponding to data $y$ and $y'$ respectively, then for every $\delta>0$, there exists $Z=Z(\delta)>0$ such that for all $y$, $y'\in\R^m$ with
    \begin{equation*}
    \max\left\{\left\|y\right\|_2,\left\|y'\right\|_2\right\}\leq\delta,
    \end{equation*}
    we have
    \begin{equation*}
    d_{Hell}(\mu^y,\mu^{y'})\leq  {Z}\left\|y-y'\right\|_{ {\Sigma_0^{-1}}}.
    \end{equation*}
    {As a result, the expectation of any polynomially bounded function $g:L_2(\Om)\to\mK$ is continuous in $y$}.
\end{enumerate}
\end{theorem}

For practical concerns, it is important to consider the finite dimensional approximation of $\mu^y$. In particular we consider the following finite dimensional approximation:
\begin{equation}\label{RNDerivativeFiniteApprox}
\f{\rd\mu_{N_1,N_2}^y}{\rd\mu_0}\propto\exp\left(-\Phi_{N_1}(u)-R_{N_2}(u)\right),
\end{equation}
where $\Phi_{N_1}(u)$ is the $N_1$ dimensional approximation of $\Phi(u)$ with $F_{N_1}$ being the $N_1$ dimensional approximation of $F$ and $R_{N_2}(u)$ is the $N_2$ dimensional approximation of $R(u)$. Theorem~\ref{Th2} provides the convergence property of $\mu_{N_1,N_2}^y$.

\begin{theorem}\label{Th2} Assume that $ {F}$ and $ {F}_{N_1}$ satisfies A1 with constants independent of $N_1$ and $R_{N_2}$ satisfies Proposition \ref{Prop1} (i)-(ii) with constants independent of $N_2$. Assume further that for all $\ep>0$, there exist two sequences $\left\{a_{N_1}(\ep)\right\}$ and $\left\{b_{N_2}(\ep)\right\}$, both of which converge to $0$, such that $\mu_0(\mX_{\ep})>0$ for all $N_1,$ $N_2\in\N$ where
\begin{equation}
\mX_{\ep}=\left\{u\in L_2(\Om):\left|\Phi(u)-\Phi_{N_1}(u)\right|\leq a_{N_1}(\ep),~\left|R(u)-R_{N_2}(u)\right|\leq b_{N_2}(\ep)\right\},
\end{equation}
then we have
\begin{equation}
d_{Hell}(\mu^y,\mu_{N_1,N_2}^y)\to0,
\end{equation}
as $N_1$, $N_2\to\infty$.
\end{theorem}
In particular, noting that $L_2(\Om)$ is a separable Hilbert space, we can consider the finite dimensional approximation of $u$, as presented in the following corollary.

\begin{corollary}\label{Coro1} Let $\left\{e_k:k\in\N\right\}$ be a complete orthonormal basis of $L_2(\Om)$. For $N\in\N$, we define
\begin{equation}
u_N=\sum_{k=1}^N\left\la u,e_k\right\ra e_k,
\end{equation}
and
\begin{equation}
\f{\rd\mu_N^y}{\rd\mu_0}\propto\exp\left(-\Phi(u_N)-R(u_N)\right).
\end{equation}
If $f$ satisfies A1-A2, then we have
\begin{equation}
d_{Hell}(\mu^y,\mu_N^y)\to0,
\end{equation}
as $N\to\infty$.
\end{corollary}

{The proof of Theorem \ref{Th2} and Corollary \ref{Coro1} can be directly followed by \cite[Theorem 2.3., Corollary 2.4]{Z.Yao2016}, hence we omit here.}

% Limited tomography introduction
\section{Application to Limited Tomography Reconstruction}
In this section, we illustrate the application of the Bayesian inference method and the NLTG prior to  limited tomography problem.
\subsection{Preliminaries}

X-ray computed tomography (CT) plays an important role in diseases diagnosis  of human body. Let $u$ denote the image to be reconstructed. Throughout this paper, we assume $u\in L_2(\mathbb{R}^2)$ is supported in a domain $\Omega$, and we only consider the two dimensional parallel beam CT for simplicity. Then the sinograme (or the projection data) $f$ is obtained by the following Radon transform~\cite{radon1917uber}:
\begin{equation}\label{RadonTransform}
  y(\alpha,s)=Pu(\alpha,s)=\int_{-\infty}^{\infty}u(s{n}+z{}^{\perp})\mathrm{d}z,~~~~~(\alpha,s)\in[0,2\pi)\times\R
\end{equation}
where ${n}=(\cos\alpha,\sin\alpha)$ and ${n}^\perp=(-\sin\alpha,\cos\alpha)$. Given a sinogram $y=Pu$ on $[0,\pi)\times\R$, the tomographic image $u$ is reconstructed via the inverse Radon transform \cite{bracewell1967inversion,Natterer1986}:
\begin{equation}\label{InverseRadonTransform}
  u(\x)=\frac{1}{2\pi^2}\int_0^\pi\int_{-\infty}^\infty\frac{1}{{x}\cdot{n}-s}\left[\frac{\partial}{\partial s}y(\alpha,s)\right]\mathrm{d}s\mathrm{d}\alpha
\end{equation}
where ${\x}=(x_1,x_2)\in\mathbb{R}^2$.

From \eqref{InverseRadonTransform}, we can easily see that the reconstruction of $u$ requires the so-called \emph{complete knowledge} of $y$ on $[0,\pi)\times\mathbb{R}$ \cite{klann2011mumford,quinto1993singularities}. However, the problem \eqref{RadonTransform} becomes ill-posed whenever the \emph{limited data} is available in the subset of $[0,\pi)\times\mathbb{R}$ due to the reduced size of detector \cite{wang2013meaning,ward2015interior} and/or the reduced number of projections \cite{pan2009commercial,E.Y.Sidky2008,K.Choi2010}. In particular, if the projection data $f$ is available on $S_r\subsetneq [0,\pi)\times \mathbb{R}$:
\begin{equation}
  S_r=\{(\alpha,s)\in[0,\pi)\times \mathbb{R}:|s|\leq r\},
\end{equation}
then there exists a nontrivial function $g$ called the ambiguity of $P$, $i.e.$ $Pg=0$ in $S_r$ \cite{wang2013meaning}. As this ambiguity $g$ is nonconstant in the region of interest (ROI) $B(0,r)$ \cite{wang2013meaning,yang2010high}, the reconstructed image via \eqref{InverseRadonTransform} using available $y$ only on $S_r$ will be deteriorated by $g$.
%With the development of the CT theories, the uniqueness and the stability of the reconstruction with limited sinogram knowledge have been completed \cite{courdurier2008solving,W.Han2009}.

In the literature, numerous  studies have been proposed to remove the ambiguity $g$ due to the restriction of $y$ on $S_r$. These studies can be classified into the following two categories: the known subregion based approaches related to the restoration of signal from the truncated Hilbert transform, and the sparsity model based approaches using the sparse approximation of tomographic images in the ROI. See \cite{T.Schuster2012,dashti2013map,heusser2012ct} for the detailed surveys. In this paper, we will use an approach based on a reference image and NLTG regularization.

\subsection{MAP and CM estimators}

In this section, we discuss how to compute the two popular point estimator in the Bayesian setting.
The first often used point estimator in the Bayesian framework is the MAP estimator,
and following the same steps as in \cite{M.Dashti2013}, we can show that the MAP estimator with the NLTG prior is the minimizer of
\begin{equation}\label{NLTGFunctional}
I(u):=\Phi(u)+\lambda J_{{NLTV}}(u)+\f{1}{2}\left\|u\right\|_{{K}}^2,
\end{equation}
over $L_2(\Om)$,
where in this problem $\Phi(u) = \frac12\|Au-y\|_{\Sigma_0}^2$
with $A$ being the linear operator in the limited tomography problem.
%Moreover it is not hard to verify that the minimizer of functional \eqref{NLTGFunctional} over $L_2(\Om)$ exists provided that ${F}$ satisfies A1-A2.
%In the present problem the mapping $F$ is the Radon transform denoted by $A$.

To minimize \eqref{NLTGFunctional}, we adopt the widely used split Bregman method \cite{goldstein2009split} which is equivalent to the alternating direction method of multipliers (ADMM). For the sake of clarity, we present the split Bregman method for \eqref{NLTGFunctional} in Algorithm \ref{alg:equivalentmap}.

\begin{algorithm}
\caption{Split Bregman Method for MAP estimate with NLTG prior.}
\label{alg:equivalentmap}
\textbf{Initialization:} $u^0 = u_{\mathrm{ref}},b^0 = d^0 =0,\ y,\ \lambda,\ \mu.$
\begin{algorithmic}[1]
   \For{$k=0,1,2,\cdots$} \\
    Update $u^{k+1}$:
    \begin{equation}\label{usubproblem}
    u^{k+1}=\mathop{\mathrm{argmin}}\limits_{u}~\Phi(u)+\frac{\mu}{2}\|d^k-\nabla_\omega u-b^k\|^2+\f{1}{2}\|u\|_K^2
    \end{equation}\\
    Update $d^{k+1}$: %*** can you check if $\|d\|_1$ in the equation below is correct?***
    \begin{equation}\label{dsubproblem}
    d^{k+1}=\mathop{\mathrm{argmin}}\limits_{d}~\lambda\|d\|_1+\frac{\mu}{2}\|d-\nabla_\omega u^{k+1}-b^k\|^2
    \end{equation}\\
    Update $b^{k+1}$:
    \begin{equation}
    b^{k+1} = b^k+ \nabla_\omega u^{k+1}-d^{k+1}.
    \end{equation}
   \EndFor
\end{algorithmic}
\end{algorithm}

To solve \eqref{usubproblem}, we use the conjugate gradient method to solve the following linear system:
\begin{equation*}
(A^T\Sigma_{0}^{-1}A-\mu \Delta_\omega +C_0^{-1})u=A^T\Sigma_{0}^{-1}y+ {\mu}\mathrm{div}_\omega(b^k-d^k).
\end{equation*}
The $d$ subproblem \eqref{dsubproblem} has the following closed form solution:
\begin{equation*}
d^{k+1}=\mathrm{shrink}(\nabla_\omega u^{k+1}+b^{k},\lambda/\mu),
\end{equation*}
where $\mathrm{shrink}(d,\lambda)$ is the shrinkage formula defined as
\begin{equation*}
\mathrm{shrink}(d,\lambda)=\max\left(\|d\|-\lambda,0\right)\f{d}{\|d\|},
\end{equation*}
with the convention that $0/0=0$.

Another often used point estimator in the Bayesian setting is the conditional mean (CM), or the posterior mean,
which is usually evaluated  using the samples drawn from the posterior distribution, often with the Markov Chain Monte Carlo (MCMC) methods.
In this work we use the preconditioned Crank-Nicolson (pCN) algorithm developed in \cite{cotter2013mcmc} for its property of being independent of
discritization dimensionality.
Simply speaking, the pCN algorithm proposes according to,
\begin{equation}
v=(1-\beta^2)^{\frac12}u+ \beta w, \label{e:pcn}
\end{equation}
where $u$ is the present position, $v$ is the proposed position, and $\beta\in[0,1]$ is the parameter controlling the stepsize,
and $w \sim \mathcal{N}(0,\Sigma_0)$.
The associated acceptance probability is
\begin{equation}
a(u,v) = \min\Bigl\{1, \exp{[\Phi(u)+J(u)-\Phi(v)-J(v)]}\Bigr\}. \label{e:acc}
\end{equation}
We describe the complete pCN algorithm in Algorithm \ref{al:pcn}.
\begin{algorithm}
 \caption{\quad The preconditioned Crank-Nicolson (pCN) algorithm}\label{al:pcn}.

    \begin{algorithmic}[1]
 \State  Initialize $u^{(0)}\in\mX $;

   \For {$n=0$ to $N-1$}

	\State Propose $v$ using Eq. ~\eqref{e:pcn};	
        \State Draw $\theta\sim U[0,1]$
        \If {$\theta\leq a(u^{(n)},v) $ }
            \State $u^{(n+1)} =  v$;
		\Else
			\State $u^{(n+1)} =  u^{(n)}$;
        \EndIf
	\EndFor
\end{algorithmic}
\medskip
\end{algorithm}

\section{Numerical Results}

\subsection{Experimental Setup}

\begin{figure}
  \centering
  \caption{XCAT images: origin, ground truth and reference with different noise levels }
  \label{Fig-1}
  \subfigure[Original image $\ u_{\mathrm{ori}}$]{\includegraphics[width=0.25\textwidth]{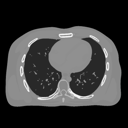}\label{fig:ori}}
  \subfigure[Ground truth $\ u_{\mathrm{gt}}$]{\includegraphics[width=0.25\textwidth]{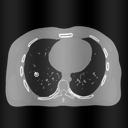}\label{fig:gt}} \\
  \subfigure[Reference image $\ u_{\mathrm{ref}}^1$]{\includegraphics[width=0.25\textwidth]{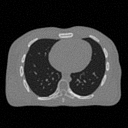}\label{fig:ref1}}
  \subfigure[Reference image $\ u_{\mathrm{ref}}^2$]{\includegraphics[width=0.25\textwidth]{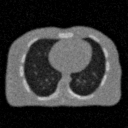}\label{fig:ref2}}
\end{figure}

In this section, we present some experimental results to demonstrate the performance of the proposed NLTG prior.
In particular we compare the results (both MAP and CM estimates) of the proposed method with those of the Filtered back projection (FBP) method,  the NLTV  regularization model and the TG method. We use the $128\times128$ XCAT image \cite{segars20104d} taking integer values in $[0,255]$ as the original image $u_{\mathrm{ori}}$. Then the ground truth image $u_{\mathrm{gt}}$ is generated by adding one round shaped object which stands for the tumor in lung and further adding the sinusoidal wave as an inhomogeneous background. Finally, we generate the reference image $u_{\mathrm{ref}}$, which can be considered as the previous CT image of the same patient taken by the same CT modality, by using $500$ projections of $u_{\mathrm{ori}}$ added by the Gaussian noise of different levels. In all experiments the reference images used are the reconstruction from the projections data with low noise level $5$ and high level $20$, which will be denoted as $u_{\mathrm{ref}}^1$ and $u_{\mathrm{ref}}^2$ respectively. Please refer to Figure \ref{Fig-1} for these images.

 {Given a reference image $u_{\mathrm{ref}}$, the covariance matrix $\Sigma_0^{-1}$ for the Gaussian measure term is computed as
\begin{equation}\label{CovMatrix}
  \Sigma_0(i,j)=\mathrm{exp}(-\frac{\|u_{\mathrm{ref}}(i)-u_{\mathrm{ref}}(j)\|^2}{h^2})
\end{equation}
following the idea of radial basis function kernel in machine learning \cite{chang2011libsvm} and the similarity weight in nonlocal means filter \cite{A.Buades2005}. We note that this choice of covariance matrix is different from usual Gaussian measure, as  the correlation  between the pixels value of $u_{\mathrm{ref}}$ is used instead of spatial distance. Such choice  aims to bring structures and edge information of the reference image to the to-be-reconstructed image, which is  especially important for reconstruction from highly missing data. In the comparison to TG method, we adopt this covariance matrix as well for a fair comparison. {As for the comparison to the NLTV method, in order to save computation and storage memory, we only use the 10 largest weights and the $4$ closest neighbors for each pixel, as adopted in \cite{X.Zhang2010}.} This will be further illustrated in the numerical results. }

To synthesize the limited projection data $y$, we generate the forward operator $A$ as the discrete Radon transform followed by the restriction onto the discretization of $[0,\pi)\times[-\frac{1}{2},\frac{1}{2}]$. Note that the size of discrete Radon transform depends on both the size of CT image and the number of projections. In this experiments, we use $50$ equally spaced projections. Then the projection data is generated by
\begin{eqnarray*}
y=Au_{\mathrm{gt}}+\eta
\end{eqnarray*}
with the Gaussian noise $\eta$ of different noise levels. Throughout this experiments, we choose the noise level to be $1$ for the low level and $10$ for the high level respectively.
Finally, the reconstructed images are further  improved by imposing the constraint of intensity  $[0, 255]$ using
\begin{equation}
u_{\mathrm{rec}}=\min\left\{\max\left(u^*,0\right),255\right\}
\end{equation}
where $u^*$ denotes the output of the four methods in comparison.
\subsection{MAP results}
%We show the prior details of $\Sigma_0^{-1}$ which plays the most important role in our NLTG model. The idea comes from the machine learning RBF kernel~\cite{chang2011libsvm},
%\begin{equation}
%  \Sigma(i,j)\ =\ \mathrm{exp}(\frac{\|u_{\mathrm{ref}}(i)-u_{\mathrm{ref}}(j)\|^2}{h^2})
%\end{equation}
%which measure the correlations between the pixel $u_{\mathrm{ref}}(i)$ and $u_{\mathrm{ref}}(j)$. According to the theories in machine learning, $\Sigma$ is a positive matrix. Hence the model (\ref{MAP_LCT}) is convex.
%\par

In solving \eqref{NLTGFunctional}, the filtering parameter $h$ in both \eqref{NonlocalWeight} and \eqref{CovMatrix} is chosen as in \cite{X.Zhang2010}. In addition, as it is not hard to see that $\Sigma_0$ is positive definite, which means that the model \eqref{NLTGFunctional} is convex, we set the maximum outer iteration number with 80 ensuring that the algorithm converges to the global minimizer. Finally, the regularization parameter $\lambda$ is manually chosen so that we can obtain the optimal restoration results.

\begin{figure}\setcounter{subfigure}{0}
  \centering
  \subfigure[Ground truth $\ u_{\mathrm{gt}}$]{\includegraphics[width=0.25\textwidth]{Groundtruth.png}}
  \subfigure[FBP]{\includegraphics[width=0.25\textwidth]{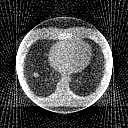}\label{fig:FBP5}}
  \subfigure[TV]{\includegraphics[width=0.25\textwidth]{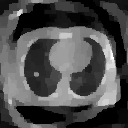}\label{fig:TV5}} \\
  \subfigure[TG with $u_{\mathrm{ref}}^1$]{\includegraphics[width=0.25\textwidth]{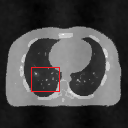}\label{TG5noisyf_ref1}}
  \subfigure[NLTV with $u_{\mathrm{ref}}^1$]{\includegraphics[width=0.25\textwidth]{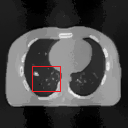}\label{NLTV5noisyf_ref1}}
  \subfigure[NLTG with $u_{\mathrm{ref}}^1$]{\includegraphics[width=0.25\textwidth]{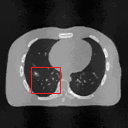}\label{NLTG5noisyf_ref1}} \\
  \subfigure[TG with $u_{\mathrm{ref}}^2$]{\includegraphics[width=0.25\textwidth]{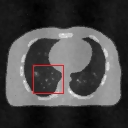}\label{TG5noisyf_ref2}}
  \subfigure[NLTV with $u_{\mathrm{ref}}^2$]{\includegraphics[width=0.25\textwidth]{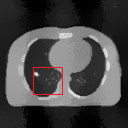}\label{NLTV5noisyf_ref2}}
  \subfigure[NLTG with $u_{\mathrm{ref}}^2$]{\includegraphics[width=0.25\textwidth]{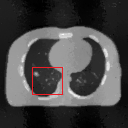}\label{NLTG5noisyf_ref2}}
  \caption{MAP results with  sinogram noise level  $5$. }
  \label{fig-2}
\end{figure}

\begin{figure}\setcounter{subfigure}{0}
  \centering
  \subfigure[Ground truth $\ u_{\mathrm{gt}}$]{\includegraphics[width=0.25\textwidth]{Groundtruth.png}\label{fig:gt}}
  \subfigure[FBP]{\includegraphics[width=0.25\textwidth]{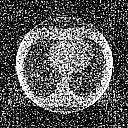}\label{fig:FBP10nosiyf}}
  \subfigure[TV]{\includegraphics[width=0.25\textwidth]{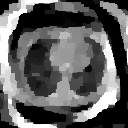}\label{fig:TV20nosiyf}} \\
  \subfigure[TG with $u_{\mathrm{ref}}^1$]{\includegraphics[width=0.25\textwidth]{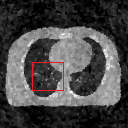}\label{TG20noisyf_ref1}}
  \subfigure[NLTV with $u_{\mathrm{ref}}^1$]{\includegraphics[width=0.25\textwidth]{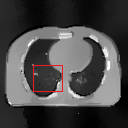}\label{NLTV20noisyf_ref1}}
  \subfigure[NLTG with $u_{\mathrm{ref}}^1$]{\includegraphics[width=0.25\textwidth]{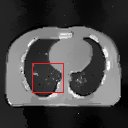}\label{NLTG20noisyf_ref1}} \\

  \subfigure[TG with $u_{\mathrm{ref}}^2$]{\includegraphics[width=0.25\textwidth]{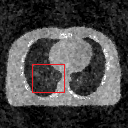}\label{NLTV20noisyf_ref2}}
  \subfigure[NLTV with $u_{\mathrm{ref}}^2$]{\includegraphics[width=0.25\textwidth]{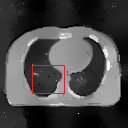}\label{NLTV20noisyf_ref2}}
  \subfigure[NLTG with $u_{\mathrm{ref}}^2$]{\includegraphics[width=0.25\textwidth]{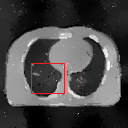}\label{NLTG20noisyf_ref2}}
  \caption{MAP results with sinogram noise level $20$.}
    \label{Fig-3}
\end{figure}

Tables \ref{table-1} and \ref{table-2} summarizes the PSNR and SSIM values of each case. As we can see from the tables, our NLTG prior consistently outperforms the other reconstruction methods, namely the FBP, TV, TG and the NLTV priors. {We present the reconstruction results with different methods in Figures \ref{fig-2} and \ref{Fig-3}}. As can be expected, both TV and FBP can not successfully reconstruct reasonable results due to limited projections. The TG prior based MAP estimate shows better performance since our proposed Gaussian covariance matrix \eqref{CovMatrix} provides some structure information from the reference image as prior.  {We can also see that, compared to the TG prior, the NLTG prior shows  the advantage of NLTV  by using the similarity in the image. In addition, compared to the NLTV prior, the NLTG prior can obtain  better recovery result, thanks to the presence of the Gaussian term which extracts more structure information by the covariance matrix computed from the reference image.}  As we can see from the figures, the visual improvements are consistent with the improvements in the indices.

\begin{table}[]
\centering
\begin{tabular}{|p{2.5cm}|p{2cm}|c|c|c|c|c|}
\hline
 Sinogram Noise Level & Reference Image & FBP  & TV & TG & NLTV & NLTG  \\ \hline
 \multirow{2}{*}{5} &$u_\mathrm{ref}^1$ &  \multirow{2}{*}{13.30}& \multirow{2}{*}{15.81} & 29.08 & 29.69 & \textbf{30.71}\\ \cline{2-2}\cline{5-7}
  &$u_\mathrm{ref}^2$ &  &   & 28.22 & 28.42 & \textbf{28.88}\\ \hline
 \multirow{2}{*}{20}&$u_\mathrm{ref}^1$ &  \multirow{2}{*}{9.40} & \multirow{2}{*}{6.46}  & 23.10 & 23.92 & \textbf{24.72}\\ \cline{2-2}\cline{5-7}
  &$u_\mathrm{ref}^2$ &    &   & 22.51 & 23.13 & \textbf{23.63}\\ \hline
\end{tabular}
\caption{MAP results: PSNR for different sinogram noise levels and reference images}
\label{table-1}
\end{table}
\begin{table}[]
\centering
\begin{tabular}{|p{2.5cm}|p{2cm}|c|c|c|c|c|}
\hline
 Sinogram Noise Level & Reference Image & FBP  & TV & TG & NLTV & NLTG  \\ \hline
 \multirow{2}{*}{5} &$u_\mathrm{ref}^1$ &  \multirow{2}{*}{0.21} & \multirow{2}{*}{0.54} & 0.88 & 0.87 & \textbf{0.91}\\ \cline{2-2}\cline{5-7}
  &$u_\mathrm{ref}^2$ &    &   & 0.87 & 0.84 & \textbf{0.86}\\ \hline
 \multirow{2}{*}{20}&$u_\mathrm{ref}^1$ &  \multirow{2}{*}{0.06} & \multirow{2}{*}{0.34} & 0.54 & 0.78 & \textbf{0.79}\\ \cline{2-2}\cline{5-7}
 &$u_\mathrm{ref}^2$ &    &   & 0.49 & \textbf{0.75} & 0.73\\ \hline
\end{tabular}
\caption{MAP results: SSIM for different sinogram noise levels and reference images}
\label{table-2}
\end{table}

%In this part, we improve the noise level in the data term to 5.0. The following tables and figures show the related results.

\subsection{CM Results}
\begin{figure}\setcounter{subfigure}{-4}
  \centering
 \subfigure{\includegraphics[width=0.24\textwidth]{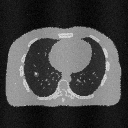}}
  \subfigure{\includegraphics[width=0.24\textwidth]{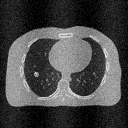}\label{f:15}}
  \subfigure{\includegraphics[width=0.24\textwidth]{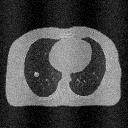}\label{f;51}}
  \subfigure{\includegraphics[width=0.24\textwidth]{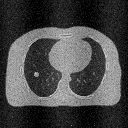}\label{f:55}} \\
   \subfigure[Noisy level of sinogram of 5 with $u_{\mathrm{ref}}^1$ ]{\includegraphics[width=0.24\textwidth]{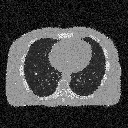}\label{f:11}}
  \subfigure[Noisy level of sinogram of 5 with $u_{\mathrm{ref}}^2$]{\includegraphics[width=0.24\textwidth]{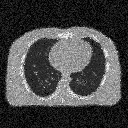}\label{f:15}}
  \subfigure[Noisy level of sinogram of 20 with $u_{\mathrm{ref}}^1$]{\includegraphics[width=0.24\textwidth]{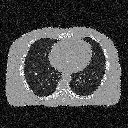}\label{f;51}}
  \subfigure[Noisy level of sinogram of 20 with  $u_{\mathrm{ref}}^2$]{\includegraphics[width=0.24\textwidth]{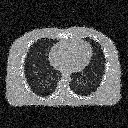}\label{f:55}}
 \caption{CM results: the condition mean for different sinogram data noise level and different references images. The upper row shows the results obtained by NLTG prior, and the down row shows  the results obtained by the TG prior.}
  \label{f:cmResultFig}
\end{figure}

The  hyperparameters in the CM model directly come from the MAP model.
For both TG prior and NLTG prior, we perform the pCN approach with $9.5\times10^5$ samples and another $0.5 \times 10^5$ samples as the pre-run. The stepwise $\beta$ has been chosen to make the acceptance probability is around $25\%$.  We show the CM results in Figure~\ref{f:cmResultFig}, and then compute the PSNR and SSIM as in Table~\ref{t:cmPSNRSSIM}.
%----------------------------------------------------------------------------------------------------------------%

 One can see from the figures and Table~\ref{t:cmPSNRSSIM}  that, the CM model with NLTG prior consistently outperformed the CM model with TG prior for different scenarios.
We can also see that the PSNR of CM model is less than that of MAP model under low sinogram noise level while the PSNR of CM model is higher than that of MAP under high sinogram noise level. Nonetheless, the results of MAP model consistently outperforms CM model in the SSIM.  {The reason of such behavior is that  MAP estimator preserves better structures and edges while CM provides smoother images that suppress the noise.  We can also see that the PSNR of CM with  different reference images under fixed sinogram noise level are almost identical, while the SSIM values are in general dependent on the reference image noise level. This suggests that, in some sense, SSIM index is more sensitive to the structures that can be brought from reference images. }

\begin{figure}\setcounter{subfigure}{-4}
  \centering
 \subfigure{\includegraphics[width=0.24\textwidth]{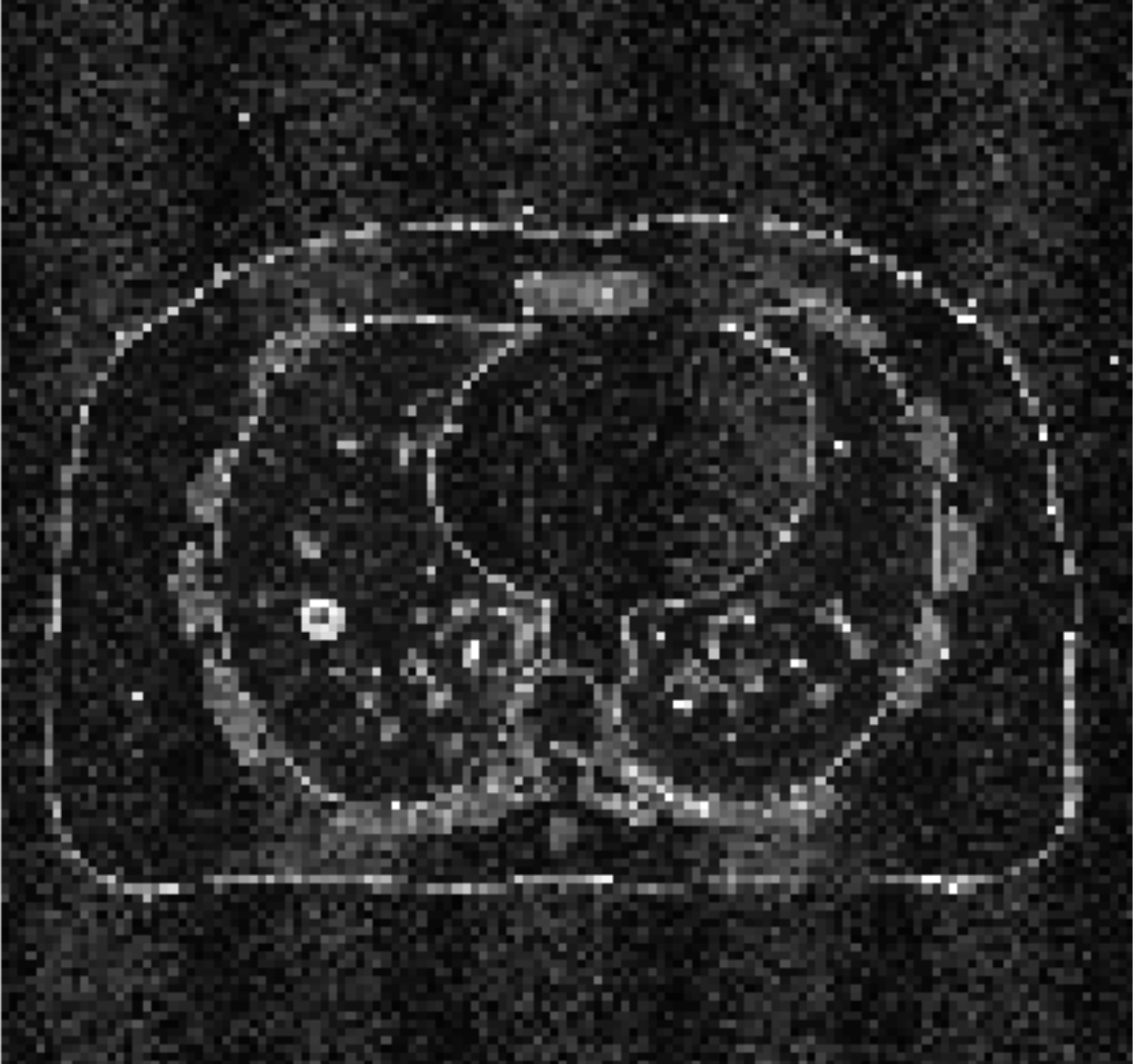}}\hspace{.01cm}
  \subfigure{\includegraphics[width=0.24\textwidth]{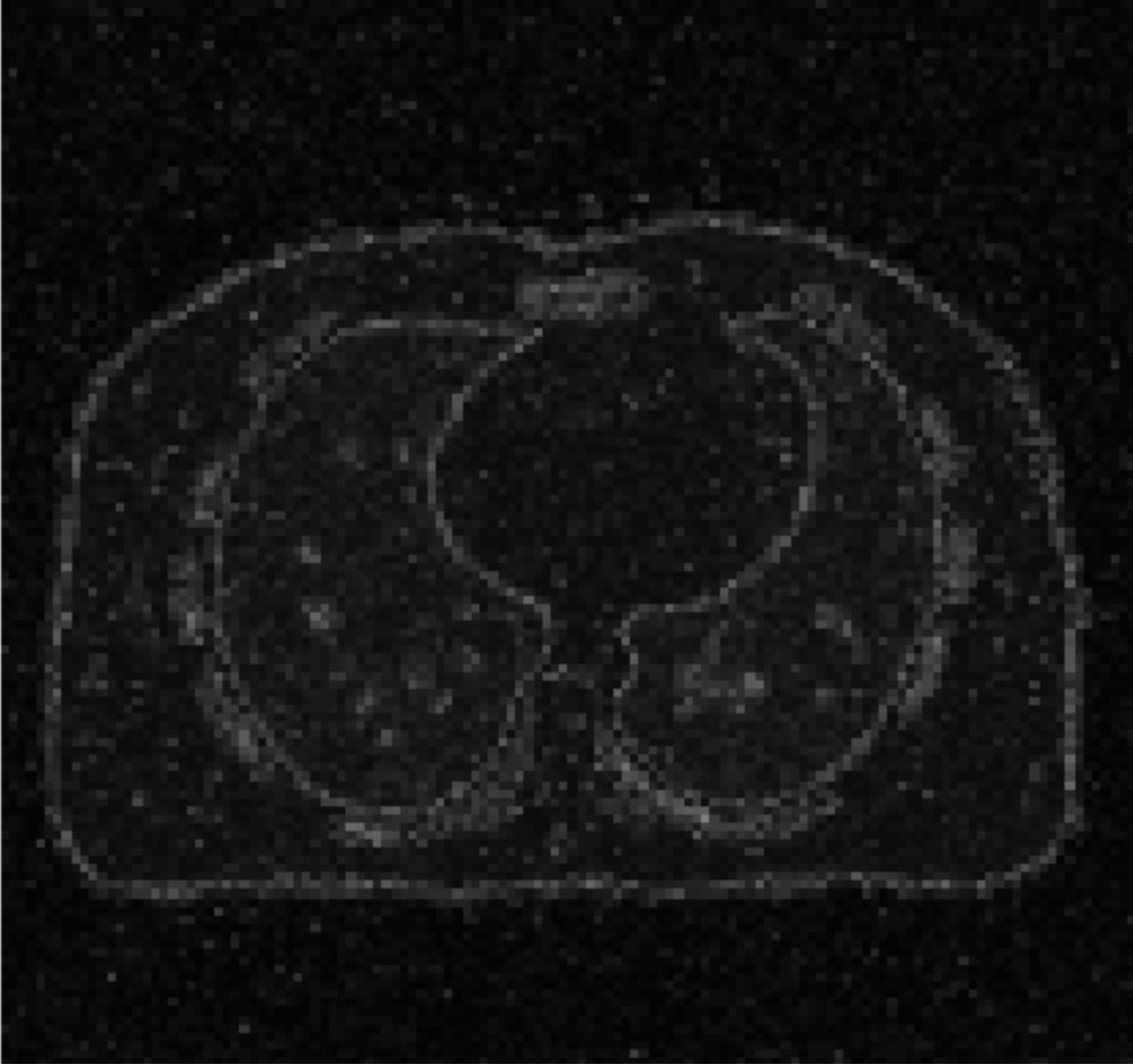}\label{f:15}}
  \subfigure{\includegraphics[width=0.24\textwidth]{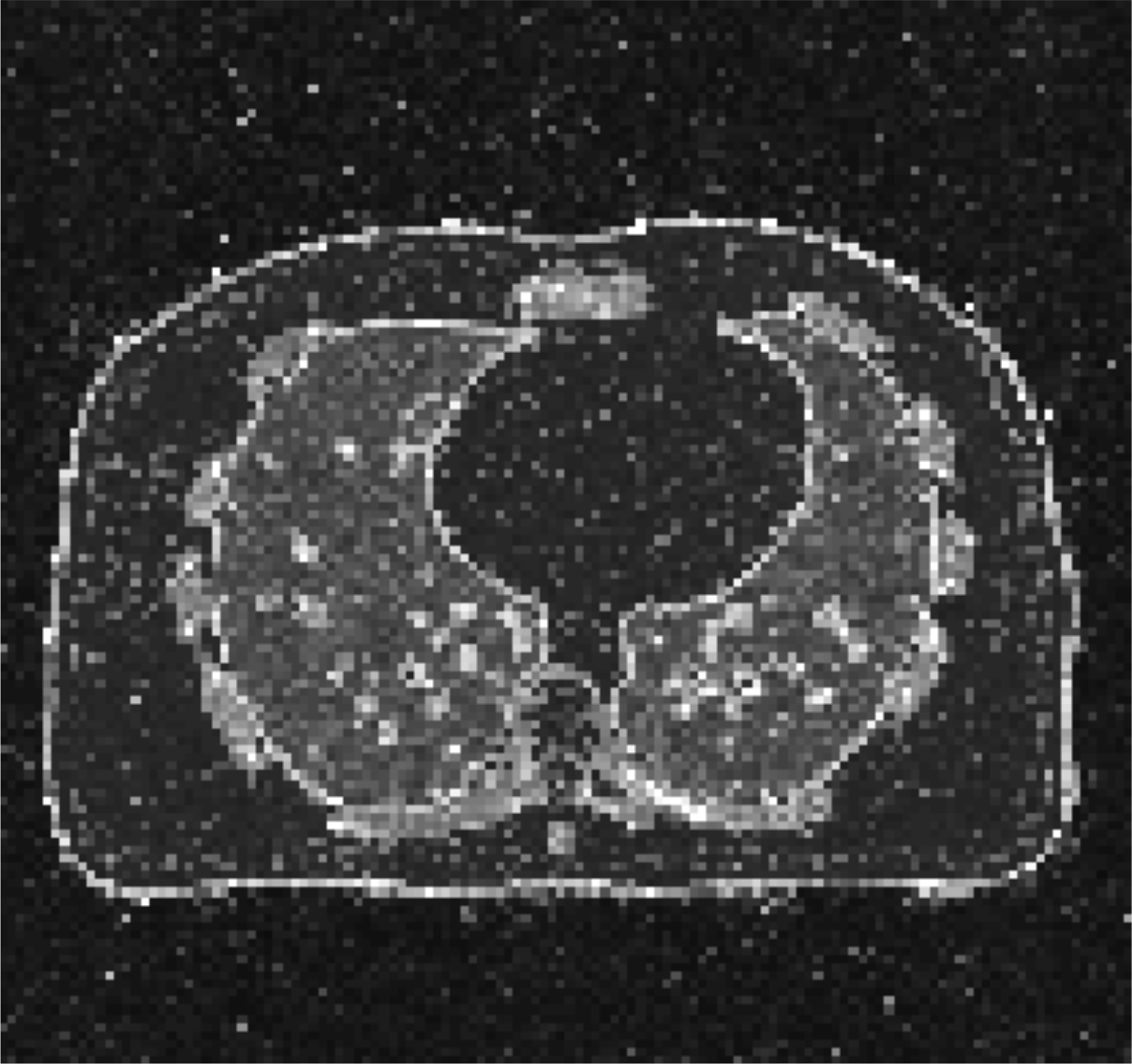}\label{f;51}}
  \subfigure{\includegraphics[width=0.24\textwidth]{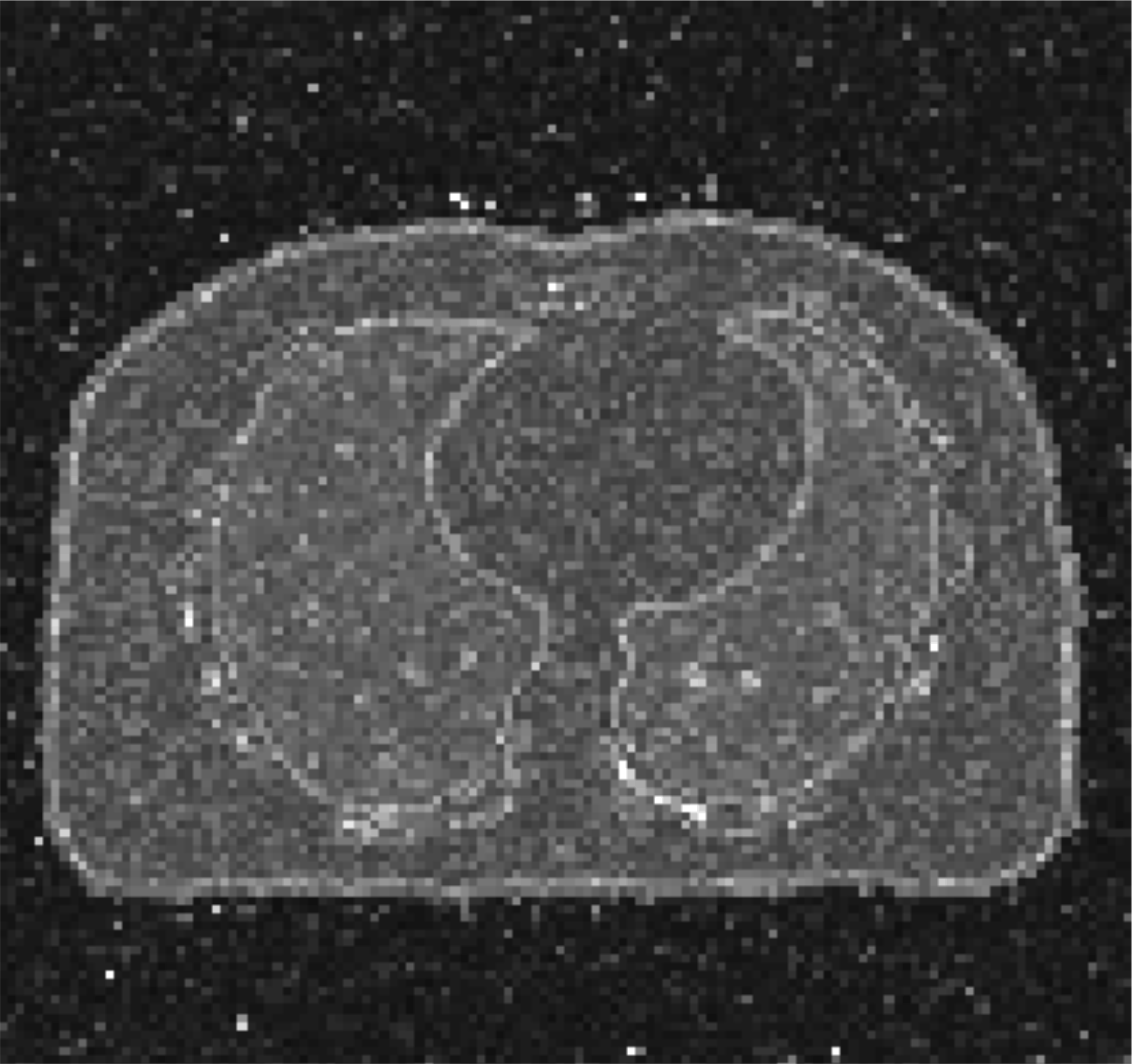}\label{f:55}} \\
   \subfigure[Noisy level of sinogram of 5 with $u_{\mathrm{ref}}^1$ ]{\includegraphics[width=0.24\textwidth]{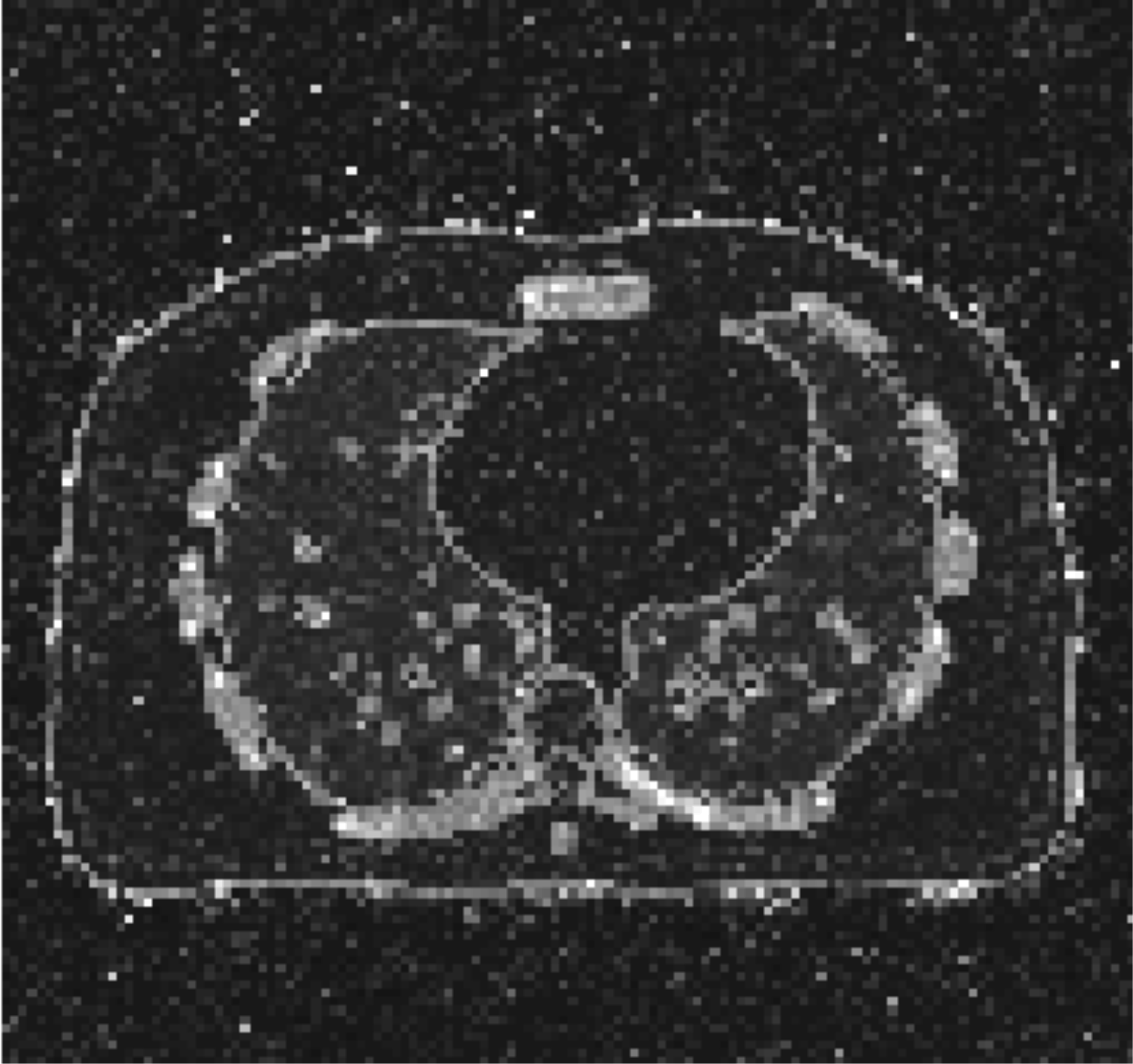}\label{f:11}}
  \subfigure[Noisy level of sinogram of 5 with $u_{\mathrm{ref}}^2$]{\includegraphics[width=0.24\textwidth]{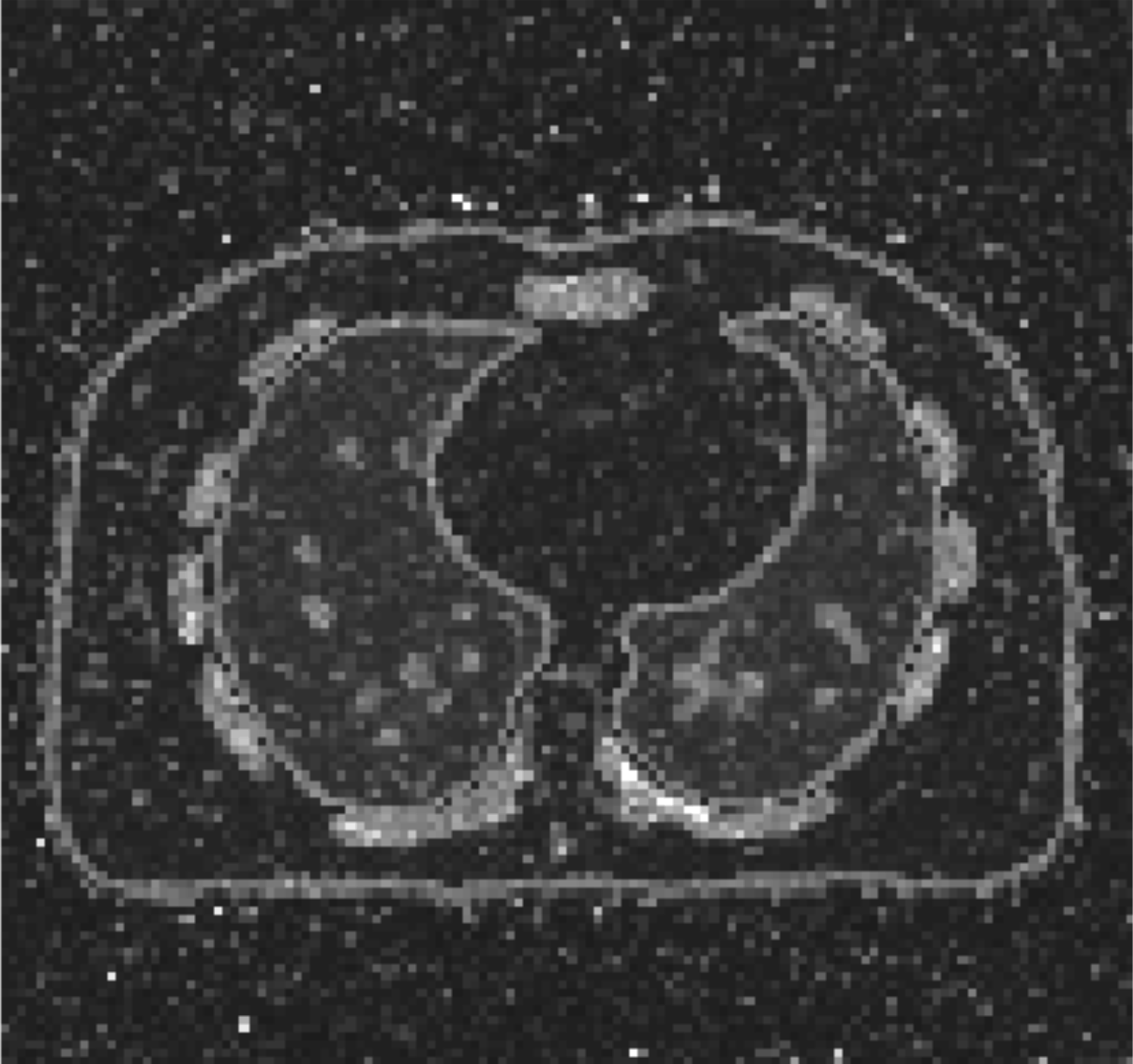}\label{f:15}}
  \subfigure[Noisy level of sinogram of 20 with $u_{\mathrm{ref}}^1$]{\includegraphics[width=0.24\textwidth]{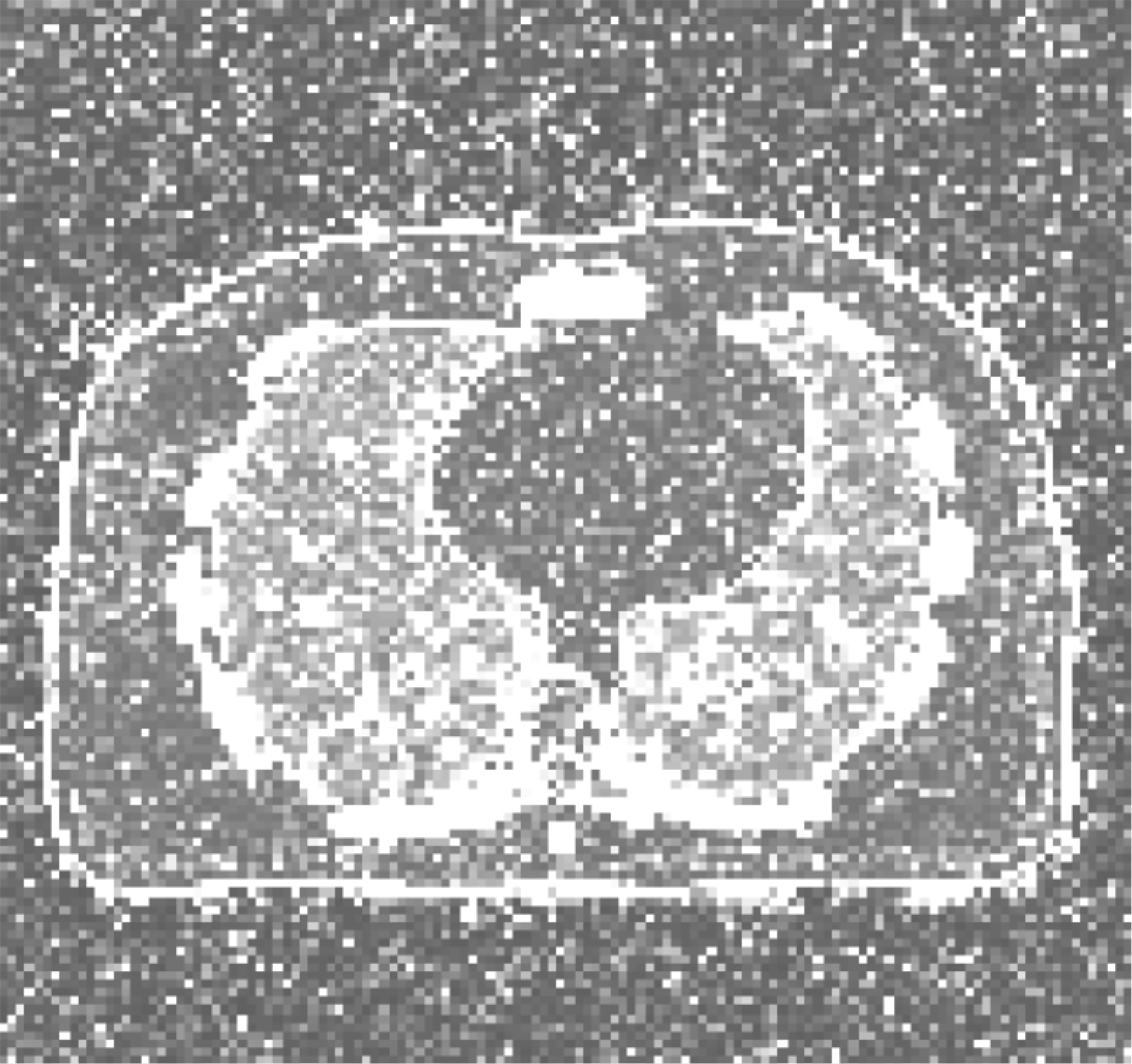}\label{f;51}}
  \subfigure[Noisy level of sinogram of 20 with  $u_{\mathrm{ref}}^2$]{\includegraphics[width=0.24\textwidth]{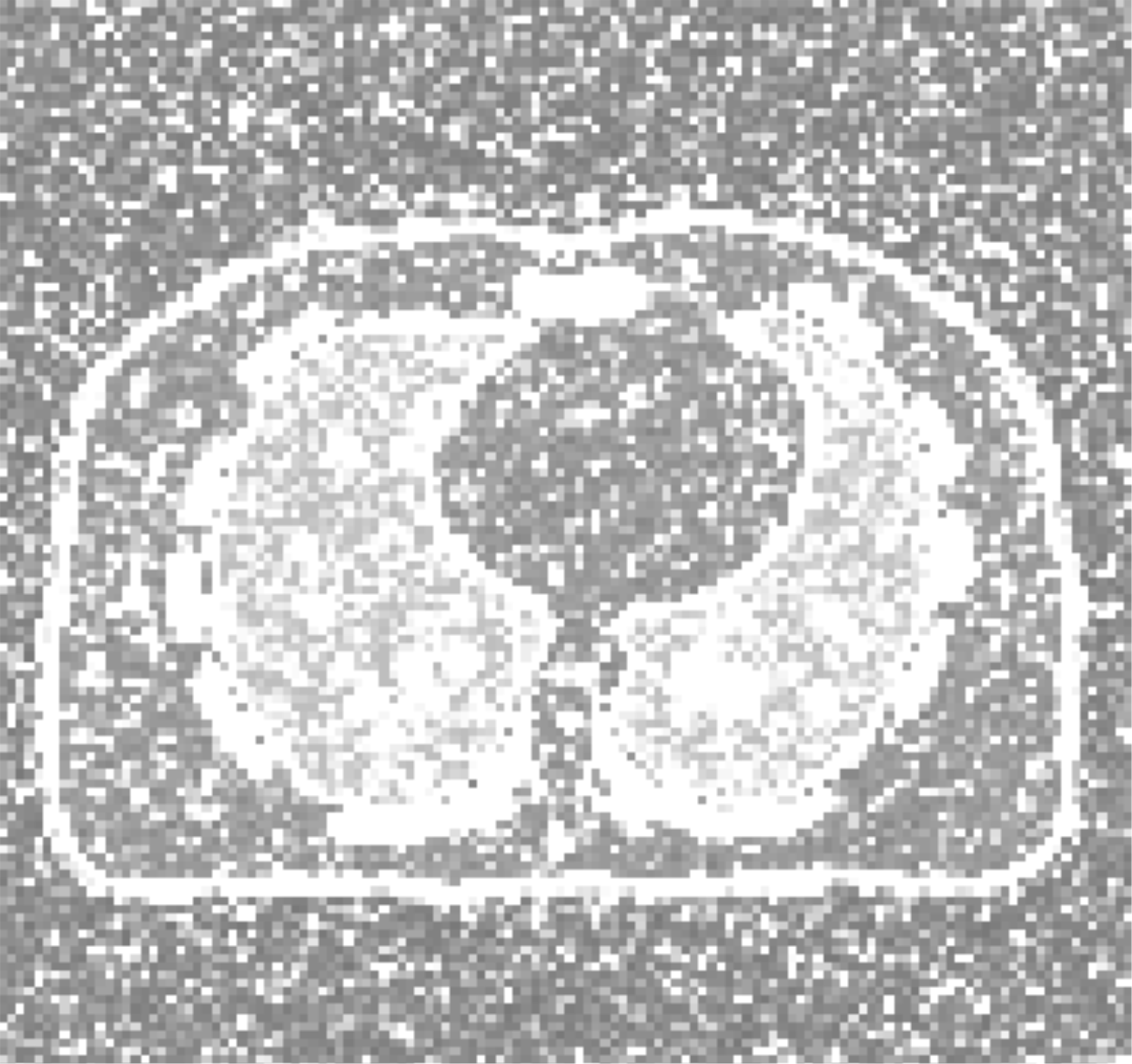}\label{f:55}}
 \caption{The 95\% confidence interval for different sinogram data noise level and different references images, the range of the values is from 0 to 100.  The upper row shows the results obtained by NLTV prior, and the down row shows  the results obtained by the TG prior. }
  \label{f:ci}
\end{figure}

The main advantage of Bayesian techniques is that it can measure the uncertainty in the estimates.  Figure~\ref{f:ci} summarizes the $95\%$ confidence interval (CI) gaps for each setting.
 Once again, the CM model with TG prior performed worse than the CM model with NLTG prior, as one would expect since the similarity of structures or distribution of an image extracted by NLTG prior can include more information than detection of local sharp edges from TG prior.

 For the  sinogram noise level of 5 ,  we can see that the CI gap of samples with reference noise level of 1 shows  higher values, especially near the edge of the image, than that of samples with reference noise level of 5.
It is worth noting that edges have high uncertainty than the smooth regions. This may be due to the fact that the NLTG prior always try to frequently tune the optimal values on the edge pixel which consequently result in high variance and  large confidence interval. Finally, we can derive a similar conclusion in the case of high sinogram noise level.

\begin{table}
\centering
\caption{CM results: SSIM and PSNR for different level of Sinogram noise and reference}
\label{t:cmPSNRSSIM}
\begin{tabular}{|p{2.5cm}|p{2cm}|c|c|c|c|}
\hline
                     &                   & \multicolumn{2}{c|}{PSNR}         & \multicolumn{2}{c|}{ SSIM}         \\
\hline
Sinogram Noise Level & Reference Image   & NLTG                           &  TG  & NLTG                            & TG   \\
\hline
\multirow{2}{*}{ 5}   & $u_\mathrm{ref}^1$ & \textbf{27.73}                          & 21.44 & \textbf{0.80}                            & 0.46  \\
\cline{2-6}
                     & $u_\mathrm{ref}^2$ & \textbf{27.90}                          & 20.95 & \textbf{0.66}                            & 0.40  \\
\hline
\multirow{2}{*}{20}  & $u_\mathrm{ref}^1$ & \textbf{25.97}                          & 19.58 & \textbf{0.62}                            & 0.37  \\
\cline{2-6}
                     & $u_\mathrm{ref}^2$ & \textbf{25.71}                          & 18.95 & \textbf{0.59}                            & 0.31  \\
\hline
\end{tabular}
\end{table}
%----------------------------------------------------------------------------------------------------------------%
\section{Conclusions}
In this paper, we consider the Bayesian inference methods for infinite dimensional inverse problems, and in particular we propose a prior distribution
that combines the nonlocal method to extract information from a reference image, with a standard Gaussian distribution. We show that the proposed prior distribution leads to a well-behaved posterior measure in the infinite dimensional setting.
We then apply the proposed method to a limited tomography problem.
The numerical experiments demonstrate the performance of the proposed NLTG prior is competitive against  existing and adapted methods,
and we also provide a comparison of the MAP and the CM estimates.
We believe that the proposed NLTG prior distribution can be useful in a large class of image reconstruction problems where reference images are available and we plan to investigate these applications in the future.

\section*{References}

\bibliographystyle{plain}

%\bibliography{ref_add_zqp}
%\bibliography{References}

\end{document}

%% file: Definition_Paper.tex
\def\R{{\mathbb R}}

\def\N{{\mathbb N}}

\def\div{\mathrm{div}}
\def\Om{\Omega}
\def\om{\omega}
\def\ep{\epsilon}

\def\f{\frac}

\def\na{\nabla}
\def\la{\langle}
\def\ra{\rangle}

\def\rd{{\mathrm d}}

\def\bi{{\mathbf i}}

\def\x{\boldsymbol{x}}
\def\y{{\boldsymbol y}}

\def\mK{{\mathcal K}}

\def\mX{{\mathcal X}}

\def\bi{\begin{itemize}} \def\ei{\end{itemize}}
\def\be{\begin{eqnarray*}}
\def\ee{\end{eqnarray*}}

\def\0{{\mathbf 0}}

\newcommand{\beq}{\begin{equation}}
\newcommand{\eeq}{\end{equation}}

\def\la{\langle}
\def\ra{\rangle}

\def\XXint#1#2#3{{\setbox0=\hbox{$#1{#2#3}{\int}$ }
\vcenter{\hbox{$#2#3$ }}\kern-.55\wd0}}

%% file: NLTG_v5_20181227.bbl
\begin{thebibliography}{10}

\bibitem{bracewell1967inversion}
RN. Bracewell and AC. Riddle.
\newblock Inversion of fan-beam scans in radio astronomy.
\newblock {\em The Astrophysical Journal}, 150:427, 1967.

\bibitem{A.Buades2005}
A.~Buades, B.~Coll, and J.~M. Morel.
\newblock {A Review of Image Denoising Algorithms, with a New One}.
\newblock {\em Multiscale Model. $\&$ Simul.}, 4(2):490--530, 2005.

\bibitem{chang2011libsvm}
C.~Chang and C.~Lin.
\newblock Libsvm: a library for support vector machines.
\newblock {\em ACM transactions on intelligent systems and technology (TIST)},
  2(3):27, 2011.

\bibitem{K.Choi2010}
K.~Choi, J.~Wang, L.~Zhu, T.~S. Suh, S.~Boyd, and L.~Xing.
\newblock {Compressed Sensing Based Cone-Beam Computed Tomography
  Reconstruction with a First-Order Method}.
\newblock {\em Med. Phys.}, 37(9):5113--5125, 2010.

\bibitem{Chung1997}
F.~R.~K. Chung.
\newblock {\em {Spectral Graph Theory}}, volume~92 of {\em CBMS Regional
  Conference Series in Mathematics}.
\newblock Published for the Conference Board of the Mathematical Sciences,
  Washington, DC; by the American Mathematical Society, Providence, RI, 1997.

\bibitem{cotter2013mcmc}
SL. Cotter, GO. Roberts, AM. Stuart, and D.~White.
\newblock {MCMC} methods for functions: modifying old algorithms to make them
  faster.
\newblock {\em Statistical Science}, 28(3):424--446, 2013.

\bibitem{M.Dashti2013}
M.~Dashti, K.~J.~H. Law, A.~M. Stuart, and J.~Voss.
\newblock {M{AP} Estimators and Their Consistency in {B}ayesian Nonparametric
  Inverse Problems}.
\newblock {\em Inverse Problems}, 29(9):095017, 27, 2013.

\bibitem{dashti2013map}
M.~Dashti, KJH. Law, AM. Stuart, and J.~Voss.
\newblock {MAP} estimators and their consistency in bayesian nonparametric
  inverse problems.
\newblock {\em Inverse Problems}, 29(9):095017, 2013.

\bibitem{Efros99}
Efros, A.~Alexei, and T.K. Leung.
\newblock Texture synthesis by non-parametric sampling.
\newblock In {\em IEEE International Conference on Computer Vision}, pages
  1033--1038, Corfu, Greece, September 1999.

\bibitem{A.Elmoataz2008}
A.~Elmoataz, O.~Lezoray, and S.~Bougleux.
\newblock {Nonlocal Discrete Regularization on Weighted Graphs: a Framework for
  Image and Manifold Processing}.
\newblock {\em IEEE Trans. Image Process.}, 17(7):1047--1060, 2008.

\bibitem{Foll1999}
G.~B. Folland.
\newblock {\em {Real Analysis: Modern Techniques and Their Applications}}.
\newblock Pure and Appl. Math. John Wiley \& Sons Inc., New York, 2nd edition,
  1999.

\bibitem{gelman2014bayesian}
A.~Gelman, JB. Carlin, HS. Stern, DB. Dunson, A.~Vehtari, and DB. Rubin.
\newblock {\em {Bayesian data analysis}}.
\newblock CRC press Boca Raton, FL, 2014.

\bibitem{G.Gilboa2007}
G.~Gilboa and S.~Osher.
\newblock {Nonlocal Linear Image Regularization and Supervised Segmentation}.
\newblock {\em Multiscale Model. $\&$ Simul.}, 6(2):595--630, 2007.

\bibitem{G.Gilboa2008}
G.~Gilboa and S.~Osher.
\newblock {Nonlocal Operators with Applications to Image Processing}.
\newblock {\em Multiscale Model. $\&$ Simul.}, 7(3):1005--1028, 2008.

\bibitem{goldstein2009split}
T.~Goldstein and S.~Osher.
\newblock The split bregman method for $l_1$-regularized problems.
\newblock {\em SIAM journal on imaging sciences}, 2(2):323--343, 2009.

\bibitem{heusser2012ct}
T.~Heu{\ss}er, M.~Brehm, S.~Marcus, S.~Sawall, and M.~Kachelrie{\ss}.
\newblock {CT} data completion based on prior scans.
\newblock In {\em Nuclear Science Symposium and Medical Imaging Conference
  (NSS/MIC), 2012 IEEE}, pages 2969--2976. IEEE, 2012.

\bibitem{kaipio2006statistical}
J.~Kaipio J and E.~Somersalo.
\newblock {\em Statistical and computational inverse problems}, volume 160.
\newblock Springer Science \& Business Media, 2006.

\bibitem{S.Kindermann2005}
S.~Kindermann, S.~Osher, and P.~W. Jones.
\newblock {Deblurring and Denoising of Images by Nonlocal Functionals}.
\newblock {\em Multiscale Model. $\&$ Simul.}, 4(4):1091--1115, 2005.

\bibitem{klann2011mumford}
E~Klann.
\newblock A mumford--shah-like method for limited data tomography with an
  application to electron tomography.
\newblock {\em SIAM Journal on Imaging Sciences}, 4(4):1029--1048, 2011.

\bibitem{lassas2004can}
M.~Lassas and S.~Siltanen.
\newblock Can one use total variation prior for edge-preserving bayesian
  inversion?
\newblock {\em Inverse Problems}, 20(5):1537, 2004.

\bibitem{li2015note}
Jinglai Li.
\newblock A note on the karhunen--lo{\`e}ve expansions for infinite-dimensional
  bayesian inverse problems.
\newblock {\em Statistics \& Probability Letters}, 106:1--4, 2015.

\bibitem{liu2015mo}
J.~Liu, H.~Ding, S.~Molloi, X.~Zhang, and H.~Gao.
\newblock Nonlocal total variation based spectral {CT} image reconstruction.
\newblock {\em Medical physics}, 42(6):3570--3570, 2015.

\bibitem{Y.Lou2010}
Y.~Lou, X.~Zhang, S.~Osher, and A.~Bertozzi.
\newblock {Image Recovery via Nonlocal Operators}.
\newblock {\em J. Sci. Comput.}, 42(2):185--197, 2010.

\bibitem{lucka2012hierarchical}
F.~Lucka, S.~Pursiainen, M.~Burger, and C.H. Wolters.
\newblock Hierarchical bayesian inference for the eeg inverse problem using
  realistic fe head models: depth localization and source separation for focal
  primary currents.
\newblock {\em Neuroimage}, 61(4):1364--1382, 2012.

\bibitem{Natterer1986}
F.~Natterer.
\newblock {\em {The Mathematics of Computerized Tomography}}.
\newblock Springer, 1986.

\bibitem{pan2009commercial}
X.~Pan, EY. Sidky, and M.~Vannier.
\newblock Why do commercial ct scanners still employ traditional, filtered
  back-projection for image reconstruction?
\newblock {\em Inverse problems}, 25(12):123009, 2009.

\bibitem{Peyre2008}
G.~Peyr\'e.
\newblock {Image Processing with Nonlocal Spectral Bases}.
\newblock {\em Multiscale Model.$\&$ Simul.}, 7(2):703--730, 2008.

\bibitem{peyre2011nonlocal}
Gabriel Peyr{\'e}, S{\'e}bastien Bougleux, and Laurent Cohen.
\newblock Non-local regularization of inverse problems.
\newblock In David Forsyth, Philip Torr, and Andrew Zisserman, editors, {\em
  Computer Vision -- ECCV 2008}. Springer Berlin Heidelberg, 2008.

\bibitem{quinto1993singularities}
ET. Quinto.
\newblock Singularities of the x-ray transform and limited data tomography in
  {R}$^{2}$ and {R}$^{3}$.
\newblock {\em SIAM Journal on Mathematical Analysis}, 24(5):1215--1225, 1993.

\bibitem{radon1917uber}
J.~Radon.
\newblock Uber die bestimmug von funktionen durch ihre integralwerte laengs
  geweisser mannigfaltigkeiten.
\newblock {\em Berichte Saechsishe Acad. Wissenschaft. Math. Phys., Klass},
  69:262, 1917.

\bibitem{rudin1992nonlinear}
L.~Rudin, S.~Osher, and E.~Fatemi.
\newblock Nonlinear total variation based noise removal algorithms.
\newblock {\em Physica D: nonlinear phenomena}, 60(1-4):259--268, 1992.

\bibitem{T.Schuster2012}
T.~Schuster, B.~Kaltenbacher, B.~Hofmann, and K.~S. Kazimierski.
\newblock {\em {Regularization Methods in {B}anach Spaces}}, volume~10 of {\em
  Radon Series on Computational and Applied Mathematics}.
\newblock Walter de Gruyter GmbH \& Co. KG, Berlin, 2012.

\bibitem{segars20104d}
W.P. Segars, G.~Sturgeon, S.~Mendonca, Jason Grimes, and Benjamin~MW Tsui.
\newblock 4d xcat phantom for multimodality imaging research.
\newblock {\em Medical physics}, 37(9):4902--4915, 2010.

\bibitem{E.Y.Sidky2008}
E.~Y. Sidky and X.~Pan.
\newblock {Image Reconstruction in Circular Cone-Beam Computed Tomography by
  Constrained, Total-Variation Minimization}.
\newblock {\em Phys. Med. Biol.}, 53(17):4777, 2008.

\bibitem{Bayesianperspective}
AM. Stuart.
\newblock {Inverse problems: a Bayesian perspective}.
\newblock {\em Acta Numer}, 19:451--559, 2010.

\bibitem{Vollmer2015}
S.~J. Vollmer.
\newblock {Dimension-Independent {MCMC} Sampling for Inverse Problems with
  Non-{G}aussian Priors}.
\newblock {\em SIAM/ASA J. Uncertain. Quantif.}, 3(1):535--561, 2015.

\bibitem{vollmer2015dimension}
Sebastian~J Vollmer.
\newblock Dimension-independent mcmc sampling for inverse problems with
  non-gaussian priors.
\newblock {\em SIAM/ASA Journal on Uncertainty Quantification}, 3(1):535--561,
  2015.

\bibitem{wang2013meaning}
G.~Wang and H.~Yu.
\newblock The meaning of interior tomography.
\newblock {\em Physics in medicine and biology}, 58(16):R161, 2013.

\bibitem{ward2015interior}
JP. Ward, M.~Lee, JC. Ye, and M.~Unser.
\newblock Interior tomography using 1d generalized total variation. part i:
  Mathematical foundation.
\newblock {\em SIAM Journal on Imaging Sciences}, 8(1):226--247, 2015.

\bibitem{yang2010high}
J.~Yang, H.~Yu, M.~Jiang, and G.~Wang.
\newblock High-order total variation minimization for interior tomography.
\newblock {\em Inverse problems}, 26(3):035013, 2010.

\bibitem{Z.Yao2016}
Z.~Yao, Z.~Hu, and J.~Li.
\newblock {A {TV}-{G}aussian Prior for Infinite-Dimensional {B}ayesian Inverse
  Problems and Its Numerical Implementations}.
\newblock {\em Inverse Problems}, 32(7):075006, 19, 2016.

\bibitem{yao2016tv}
Zhewei Yao, Zixi Hu, and Jinglai Li.
\newblock A tv-gaussian prior for infinite-dimensional bayesian inverse
  problems and its numerical implementations.
\newblock {\em Inverse Problems}, 32(7):075006, 2016.

\bibitem{X.Zhang2010}
X.~Zhang, M.~Burger, X.~Bresson, and S.~Osher.
\newblock {Bregmanized Nonlocal Regularization for Deconvolution and Sparse
  Reconstruction}.
\newblock {\em SIAM journal on imaging sciences}, 3(3):253--276, 2010.

\bibitem{X.Zhang2010a}
X.~Zhang and T.~F. Chan.
\newblock {Wavelet Inpainting by Nonlocal Toral Variation}.
\newblock {\em Inverse problems and Imaging}, 4(1):191--210, 2010.

\bibitem{D.Zhou2005}
D.~Zhou and B.~Sch{\"o}lkopf.
\newblock {\em {Regularization on Discrete Spaces}}.
\newblock Springer Berlin Heidelberg, Berlin, Heidelberg, 2005.

\end{thebibliography}
